\documentclass[12pt]{article}
\usepackage{array,amsmath,amssymb}
\usepackage{amsthm}
\usepackage{fullpage}
\usepackage{url}

\newtheorem{theorem}{Theorem}
\newtheorem{lemma}[theorem]{Lemma}
\newtheorem{corollary}[theorem]{Corollary}

\newtheorem{question}[theorem]{Question}
\newtheorem{example}[theorem]{Example}
\newtheorem{proposition}[theorem]{Proposition}

\newcommand{\peb}[1]{\Pi\left(#1\right)}
\newcommand{\optpeb}[1]{\Pi_{OPT}\left(#1\right)}
\newcommand{\optpebf}{\Pi_{OPT}}

\newenvironment{Proof}{\par\noindent{\bf Proof.} }{\qed\par\bigskip}

\input epsf.tex
\def\qed{\hfill{\setlength{\fboxsep}{0pt}\framebox[7pt]{\rule{0pt}{7pt}}}\newline}
\def\st{\colon\,}   
\def\NN{{\Bbb N}}
\def\cP{{\cal L}}
\def\bG{{\bold G}}
\def\({\left(}  \def\){\right)}

\def\FR#1#2{\frac{#1}{#2}}
\def\FL#1{\left\lfloor{#1}\right\rfloor}
\def\CL#1{\left\lceil{#1}\right\rceil}
\def\floor#1{\FL{#1}}
\def\ceil#1{\CL{#1}}
\def\MAP#1#2#3{#1\colon\,#2\to#3}
\def\esub{\subseteq}
\def\VEC#1#2#3{#1_{#2},\ldots,#1_{#3}}
\def\SE#1#2#3{\sum_{#1=#2}^{#3}}
\def\cart{\>\hbox{${\vcenter{\vbox{
   \hrule height 0.4pt\hbox{\vrule width 0.4pt height 4.5pt
\kern4pt\vrule width 0.4pt}\hrule height 0.4pt}}}$}\>}

\begin{document}

\title{Pebbling and Optimal Pebbling in Graphs}
\author{David P. Bunde \thanks{Computer Science Department,
University of Illinois, Urbana, IL 61801, bunde@uiuc.edu.
Partially supported by NSF grant CCR 0093348.}
\and Erin W. Chambers\thanks{Computer Science Department,
University of Illinois, Urbana, IL 61801, erinwolf@uiuc.edu.
Supported by a National Science Foundation Graduate Research Fellowship.}
\and Daniel Cranston\thanks{Computer Science Department,
University of Illinois, Urbana, IL 61801, dcransto@uiuc.edu.}
\and Kevin Milans\thanks{Computer Science Department,
University of Illinois, Urbana, IL 61801, milans@uiuc.edu.}
\and
Douglas B. West\thanks{Mathematics Department,
University of Illinois, Urbana, IL 61801, west@math.uiuc.edu.
Work supported in part by the NSA under Award No.~MDA904-03-1-0037.}
}

\date{}
\maketitle

\begin{abstract}
Given a distribution of pebbles on the vertices of a graph $G$, a
{\it pebbling move} takes two pebbles from one vertex and puts one on
a neighboring vertex.  The {\it pebbling number} $\peb{G}$ is the minimum $k$
such that for every distribution of $k$ pebbles and every vertex $r$, it is
possible to move a pebble to $r$.  The {\it optimal pebbling number}
$\optpeb{G}$ is the minimum $k$ such that some distribution of $k$ pebbles
permits reaching each vertex.

We give short proofs of prior results on these parameters for paths, cycles,
trees, and hypercubes, a new linear-time algorithm for computing $\peb{G}$ on
trees, and new results on $\optpeb{G}$.  If $G$ is a connected $n$-vertex
graph, then $\optpeb{G}\le\CL{2n/3}$, with equality for paths and cycles.
If $\bG$ is the family of $n$-vertex connected graphs with minimum degree $k$,
then $2.4\le \max_{G\in \bG} \optpeb{G} \FR{k+1}n\le 4$ when $k>15$ and $k$ is
a multiple of 3.  Finally, $\optpeb{G}\le 4^tn/((k-1)^t+4^t)$ when $G$ is a
connected $n$-vertex graph with minimum degree $k$ and girth at least $2t+1$.
For $t=2$, a more precise version of this last bound is
$\optpeb{G}\le 16n/(k^2+17)$.
\end{abstract}

\section{Introduction}

Graph pebbling is a model for the transmission of consumable resources.
Initially, pebbles are placed on the vertices of a graph $G$ according to a
{\it distribution} $D$, a function $\MAP D{V(G)}{\NN\,\cup\{0\}}$.
A {\em pebbling move} from a vertex $v$ to a neighbor $u$ takes away two
pebbles at $v$ and adds one pebble at $u$.  Before the move, $v$ must have at
least two pebbles.  A {\em pebbling sequence} is a sequence of pebbling moves.

Given a distribution and a ``root'' vertex $r$, the task is to put a pebble on
$r$.  A distribution $D$ is {\em $r$-solvable} (and $r$ is {\it reachable}
under $D$) if $r$ has a pebble after some (possibly empty) pebbling sequence
starting from $D$.  For a graph $G$, let $\peb{G,r}$ be the least $k$
such that every distribution of $k$ pebbles on $G$ is $r$-solvable.  A
distribution $D$ is {\em solvable} if every vertex is reachable under $D$.  The
{\em pebbling number} of a graph $G$, denoted $\peb{G}$, is the least $k$ such
that every distribution of $k$ pebbles on $G$ is solvable.  The
{\em optimal pebbling number} of $G$, denoted $\optpeb{G}$, is the least $k$
such that some distribution of $k$ pebbles is solvable.

Graph pebbling originated in efforts of Lagarias and Saks to shorten a result
in number theory.  A survey by Hurlbert~\cite{hurlbert99} describes this
history and summarizes early results.  Hurlbert introduced a useful
generalization.  A distribution $D$ is {\em $m$-fold $r$-solvable} (and $r$ is
{\it $m$-reachable} under $D$) if $r$ has at least $m$ pebbles after some
(possibly empty) pebbling sequence.  A distribution $D$ is $m$-fold
solvable if every vertex is $m$-reachable under $D$.  When $m=2$, we say that
an $m$-fold solvable distribution is {\it doubly solvable}.

Moews~\cite{moews92} developed several useful tools for computing pebbling
numbers.  (An unpublished longer version of the paper \cite{moews92} appears on
his webpage \cite{moews}.)
We call the first of these tools the {\em Weight Argument}, which we express
here for $m$-fold solvability.  Given a root $r$ and distribution $D$, let
$a_{i,r}$ be the total number of pebbles on vertices at distance $i$ from $r$.
A pebbling move cannot increase the sum $\sum_{i\geq 0}a_{i,r}2^{-i}$.
Therefore, $m$-fold $r$-solvability of $D$ requires the {\em weight inequality}
$\sum_{i\geq0}a_{i,r}2^{-i}\geq m$.

Our other main tool is that when each pebbling move is represented by a
directed edge from the vertex losing pebbles to the vertex gaining a pebble,
no directed cycle is needed.  If $r$ is reachable using moves containing a
cycle, then also $r$ is reachable using a proper subset of these moves.  In
particular, if a distribution is $r$-solvable, then $r$ is reachable without
moving a pebble in both directions along any edge.

To make this precise, say that a directed multigraph $H$ is {\it orderable}
under a distribution $D$ if some linear ordering $\sigma$ of $E(H)$ is a valid
pebbling sequence starting from $D$.  For such $D$ and $H$, the {\it balance}
of a vertex $v$ is $d_H^-(v)+D(v)-2d_H^+(v)$, where $d_H^-(v)$ and $d_H^+(v)$
are the indegree and outdegree of $v$ under $H$.  When $H$ is orderable under
$D$ (by $\sigma$), each vertex has nonnegative balance, since the balance is the
number of pebbles at $v$ after applying $\sigma$.  The {\em No-Cycle Lemma}
states that if $H$ is orderable under $D$, then it has an acyclic subgraph $H'$
such that $H'$ is orderable under $D$ and gives balance to each vertex at least
as large as does $H$.  The lemma was proved in \cite{nocycle} and earlier in
\cite{moews92} and has a short proof in \cite{milans}.

The pebbling number is known exactly for some special graphs.
Moews~\cite{moews92} observed that a distribution on a path rooted at its end
is solvable if and only if the weight inequality holds; thus
$\peb{P_n}=2^{n-1}$ for the $n$-vertex path $P_n$, since each pebble
contributes weight at least $2^{-(n-1)}$.  The $n$-vertex cycle $C_n$ is more
complicated; Pachter et al~\cite{pachter95} proved that $\peb{C_{2k}}=2^k$ and
$\peb{C_{2k+1}}=2\floor{2^{k+1}/3}+1$.  For the $k$-dimensional hypercube
$Q_k$, Chung~\cite{chung89} proved that $\peb{Q_k}=2^k$.  For a rooted
tree, Moews~\cite{moews92} showed how to calculate the pebbling number from
decompositions into paths.  For general graphs, Milans and Clark~\cite{milans}
showed that recognizing $\peb{G}\le k$ is a $\Pi_2^{\rm P}$-complete problem,
meaning that it is complete for the class of languages computable in polynomial
time by coNP machines equipped with an oracle for an NP-complete language.

Study of the optimal pebbling number began with the result of Pachter
et al~\cite{pachter95} that $\optpeb{P_n}=\ceil{2n/3}$.  Moews~\cite{moews98}
proved that $(4/3)^k\le\optpeb{Q_k}\le(4/3)^{k+O(\log k)}$ and proved a related
result for $\optpebf$ on cartesian product graphs.  Milans and
Clark~\cite{milans} proved that computing $\optpebf$ is NP-hard on
arbitrary graphs.

In this paper, we present several new results (mostly on optimal pebbling)
and simpler proofs for some previously-known results.  The proofs are simple
due to lemmas that restrict the form of pebble distributions that need to be
considered.  These lemmas and another lower-bound technique for optimal
pebbling comprise our main contribution.  Our undirected graphs are simple and
connected.  We use $V(G)$ and $E(G)$ to denote the vertex set and edge set of a
graph $G$, with sizes $n(G)$ and $e(G)$.

For the pebbling number, we give another proof of the result of
Moews~\cite{moews92} on calculating $\peb{T,r}$ from a particular decomposition
of a tree $T$ into paths.  We extend his result to give a linear-time algorithm
for computing $\peb{T}$.  Also, we give short proofs of the results of Pachter
et al~\cite{pachter95} that $\peb{C_{2k}}=2^k$ and
$\peb{C_{2k+1}}=2\floor{2^{k+1}/3}+1$.

Our approach (especially for paths and cycles), relies on a precise version of
the following intuition.  For distributions with $k$ pebbles, the hardest ones
to make solvable are concentrated on one or two vertices, while the easiest
ones are spread over many vertices.  Thus to determine $\peb{G}$ we consider
concentrated distributions, while to determine $\optpeb{G}$ we consider
``smooth'' distributions.

For optimal pebbling, the {\it Smoothing Lemma} is that for each graph a
solvable distribution of minimum size exists with at most two pebbles on each
vertex of degree at most 2.  This leads to a simpler proof of the result of
Pachter et al~\cite{pachter95} that $\optpeb{P_n}=\ceil{2n/3}$ and a proof that
$\optpeb{C_n}=\ceil{2n/3}$.  (We recently learned that Friedman and
Wyels~\cite{friedman} have obtained another short derivation of $\optpeb{P_n}$
different from ours, and like us they adapted those ideas to compute
$\optpeb{C_n}$.)

We also show that $\optpeb{T} \leq \ceil{2n/3}$ for every $n$-vertex tree $T$,
which immediately yields $\optpeb{G} \leq \ceil{2n/3}$ for every connected
$n$-vertex graph $G$, and we give a short proof of the result of
Moews~\cite{moews98} that $\optpeb{Q_k} \geq (4/3)^k$.

Let $G$ be a connected $n$-vertex graph with minimum vertex degree $k$.
Czygrinow~\cite{czygrinow} observed that $\optpeb{G}\le 4\FR{n}{k+1}$.  We
construct families of such graphs with $\optpeb{G}\ge 2\FR{n}{k+1}$ for all $k$
and with $\optpeb{G}\ge(2.4-\frac{24}{15k+5})\FR{n}{k+1}$ when $k$ is divisible
by 3.  These results use another lower-bound technique, the simplest version of
which is that if $G$ is obtained from $H$ by collapsing sets of vertices into
single vertices, then $\optpeb{H}\ge \optpeb{G}$.

We obtain tighter bounds when we further restrict $G$ to have girth (minimum
cycle length) at least $2t+1$.  Suppose that $k\ge3$ and $t\ge2$, but exclude
the case $(k,t)=(3,2)$.  Letting $c_k(t)=1+k\SE i1t(k-1)^{i-1}$ and
$c'(t)=(4^t-2^{t+1})\frac t{t-1}$, we prove that
$\optpeb{G}\le 4^tn/(c_k(t)+c'(t))$.  When $G$ has girth at
least 5, this yields $\optpeb{G}\le \FR{16n}{k^2+17}$.
Among graphs with girth 4, we show that
$\optpeb{C_m\cart K_2}=\optpeb{P_m\cart K_2}= m$ when $m\ge3$ (except $m+1$
when $m=5$), where $\cart$ denotes cartesian product (see
Section~\ref{girth-sec}).  The same bound holds also for the graph consisting
of a $2m$-cycle with chords added joining opposite vertices.

Our results on pebbling number of trees and pebbling number of cycles
appear in Sections~\ref{computing-pebnum-section} and \ref{pebcycles},
respectively.  The final three sections discuss optimal pebbling number and
present the new results.  In addition to the results mentioned above, we pose
the question of whether every connected $n$-vertex graph with minimum degree at
least 3 has optimal pebbling number at most $\CL{n/2}$.

\section{Pebbling Number of Trees} \label{computing-pebnum-section}

Moews~\cite{moews92} showed how to compute the pebbling number of a tree from a
decomposition into paths.  In this section, we prove this more simply and show
how to find an optimal decomposition in linear time.

A partition of the edge set of a tree is a {\em path partition} if each set
in the partition is a (directed) path when all edges are directed toward a root
$r$.  The {\em length list} of a path partition is the list of path lengths in
nonincreasing order.  A path partition {\em majorizes} another if its length
list is larger than the other's in the first position where they differ.
Majorization is a linear (lexicographic) order on length lists, but distinct
path partitions may have length lists that are the same.  A path partition with
root $r$ is {\em $r$-optimal} if it is not majorized by any other path
partition with root $r$.  It is {\em optimal} if it is not majorized by any
path partition with any root.  We use {\it leaf} to describe a vertex of degree
1 in any graph.

Moews~\cite{moews92} showed how to determine $\peb{T,r}$ from an optimal path
partition of a tree $T$ rooted at a vertex $r$.  Our proof is shorter and
simpler.

\begin{theorem}[{\rm Moews~\cite{moews92}}] \label{peb-num-tree-theorem}
If the length list of an $r$-optimal path partition of tree $T$ with root $r$
is $l_1, \ldots, l_m$, then
\[
\peb{T,r} = \(\sum_{i=1}^m 2^{l_i}\) - m + 1.
\]
\end{theorem}

\begin{Proof}   
We have observed $r$-solvability never requires moving pebbles in both
directions along an edge.  Thus in a tree we may direct all edges toward
the root and move pebbles only in that direction.  Let $\cP$ be an optimal path
partition of a tree $T$ rooted at $r$, and let $(l_1,\ldots,l_m)$ be the length
list of $\cP$.

{\it Lower Bound.}
We construct a non-$r$-solvable distribution with
$\sum_{i=1}^m\left(2^{l_i}-1\right)$ pebbles.
If some path in $\cP$ starts at a nonleaf vertex, then another path ends there,
and they combine to produce a path partition majorizing $\cP$.  Hence in $\cP$
each path begins at a leaf.  For each path of length $l_i$ in $\cP$, we put
$2^{l_i}-1$ pebbles on the starting leaf.  Now no pebble can be the first
pebble to reach the end of the path in $\cP$ on which it starts.  Hence pebbles
never reach the end or move off their starting path in $\cP$.  In particular,
no pebble can reach $r$.

{\it Upper Bound.}
We show that every distribution with more than $\(\sum_{i=1}^m2^{l_i}\)-m$
pebbles is $r$-solvable, using a weight function based on $\cP$.  Let $P_i$ be
the path in $\cP$ corresponding to length $l_i$.  Given a distribution $D$, let
$a_{i,j}$ be the number of pebbles on $P_i$ at distance $j$ from the end.
Let $w_i(D)=2^{l_i}\sum_{j=1}^{l_i}a_{i,j}2^{-j}$, and let
$w(D)=\sum_{i=1}^m w_i(D)$.

The function $w_i(D)$ differs from the standard weight function on a path in
two ways: we multiply by an extra factor of $2^{l_i}$, and we sum over $j\ge1$
rather than $j\ge0$.  We sum over $j\geq 1$ because the end of $P_i$ is inside
a longer path (or is $r$); we avoid counting pebbles twice.  The factor of
$2^{l_i}$ ensures that moves toward $r$ do not decrease the total weight.

A path $P_i$ in $\cP$ is {\it full} under distribution $D$ if
$w_i(D)\ge2^{l_i}$.  If $P_i$ is full, then we can move a pebble along $P_i$ to
its end, where it will be on $r$ or contribute to the weight of a path that
extends closer to $r$.  Each move within a path does not change the total
weight.  When a pebble moves from $P_i$ to $P_{i'}$, the weight decreases by
$2\cdot 2^{l_i-1}$ and increases by $2^{l_{i'}-j}$, where $j$ is the distance
from the new location to the end of $P_{i'}$.  If $l_i>l_{i'}-j$, then $P_i$ can
replace the beginning of $P_{i'}$ to produce a path partition majorizing $\cP$;
the optimality of $\cP$ prevents this.  Hence $l_i\le l_i'-j$, and moving a
pebble from $P_i$ to $P_{i'}$ does not decrease the weight.

Given an optimal path partition with lengths $\VEC l1m$, let $D$ be a
distribution under which $r$ is not reachable.  If $D$ has more than
$\sum_{i=1}^m\(2^{l_i}-1\)$ pebbles, then by the pigeonhole principle some path
is full, since each pebble not on $r$ contributes at least 1 to the weight of
the tree.  We have shown that no move toward $r$ decreases the total weight,
except when a pebble is moved onto $r$ and no longer contributes.  Every
pebbling sequence terminates, since each move reduces the total number of
pebbles.  Since the weight never decreases, the sequence can only terminate by
moving a pebble onto $r$.
\end{Proof}

In his survey \cite{hurlbert99}, Hurlbert attributes the corollary to Moews.

\begin{corollary}\label{treecor}
If the length list of an optimal path partition of tree $T$ is
$l_1, \ldots, l_m$ then
\[
\peb{T} = \sum_{i=1}^m 2^{l_i} - m + 1.
\]
\end{corollary}
\begin{Proof}
Since exponentiation is a convex function, the formula in
Theorem~\ref{peb-num-tree-theorem} is maximized by an $r$-optimal path
partition.  Also $\peb{T}=\max_{r\in V(T)}\peb{T,r}$.
Hence the claim follows from Theorem~\ref{peb-num-tree-theorem}.
\end{Proof}

The difficulty in applying Corollary~\ref{treecor} is in finding
an optimal path partition.  Given a root, a natural idea is to select a longest
path greedily and iterate.  Although this works, it disconnects the tree,
leaving awkward bookkeeping details.  The inductive proof is simpler if we peel
away shorter paths first.  A {\it peripheral vertex} in a tree is an endpoint
of a longest path.  A {\it branch vertex} in a tree is a vertex of degree at
least 3.  An {\it $x,y$-path} in a graph is a path with endpoints $x$ and $y$.

\begin{theorem}\label{tree}
There is a linear-time algorithm to compute the pebbling number of trees.
In particular, if $r$ is an endpoint of a longest path in $T$, then
$\peb{T,r}=\peb{T}$, and any longest path to $r$ can be chosen as a
path in an $r$-optimal path partition.
\end{theorem}
\begin{Proof}
In a tree, the vertices at greatest distance from a vertex $x$ are endpoints of
a longest path.  Hence a single breadth-first search from an arbitrary vertex
finds a peripheral vertex $r$.  Another breadth-first search from $r$ finds
a longest path $R$, ending at another vertex $r'$.

With $R$ chosen, another breadth-first search computes distances from $R$.
We find an $r$-optimal path partition using these distances.  The partition
will have $R$ as a path, and it will be both $r$-optimal and $r'$-optimal.
We view all edges off $R$ as directed toward $R$.

Suppose that $R$ is not all of $T$.  Iteratively, we select a leaf $x$ closest
to $R$ among the leaves that remain in the tree.  Let $y$ be the closest branch
vertex to $x$ in $T$; vertex $y$ is well-defined.  Since $R$ is a longest path,
$y$ cannot be $r$ or $r'$.  Let $P$ be the $x,y$-path in $T$.  Put $P$ into the
path partition and delete $P$ from the tree, leaving only the endpoint $y$.
When the remaining tree is just $R$, it becomes the last path in the partition.
(We can pause the computation of distances from $R$ each time a leaf is found
and extract $P$ then.)

We prove, by induction on the number of vertices outside $R$, that the path
$P$ deleted at each step lies in an $r$-optimal path partition of the tree
remaining at that step.  By the majorization criterion, the path $P'$
containing $x$ in an $r$-optimal path partition $\cP$ contains all of $P$.
If $P'$ continues past $y$, then some path $Q$ in $\cP$ ends at $y$.  We have
observed that $Q$ starts at a leaf, so $Q$ is at least as long as $P$,
by the choice of $x$.

Let $Q'$ be the union of $Q$ and the part of $P'$ after $y$.  Let $\cP'$ be the
partition obtained from $\cP$ by replacing $P'$ and $Q$ with $P$ and $Q'$.  Now
$P$ and $Q'$ are shortest and longest, respectively, among $\{P,P',Q,Q'\}$.
If $Q$ is longer than $P$, then $\cP'$ majorizes $\cP$.  Otherwise, $\cP'$ and
$\cP$ have the same length list; hence $\cP$ is an $r$-optimal path partition
containing $P$.

Thus $P$ occurs in an $r$-optimal path partition $\cP$.  The remainder of
$\cP$ is an $r$-optimal path partition of the remaining tree $T'$.  Distances
from $R$ are the same in $T'$ as in $T$.  By the induction hypothesis, the
remainder of the algorithm produces an $r$-optimal path partition of $T'$ that
contains $R$.  It combines with $P$ to yield the desired path partition of $T$.

The partition we have produced is also $r'$-optimal, since the computation is
the same when viewed from $r'$ (distances from $R$ are the same).

Since we have found an $r$-optimal path partition containing a longest path,
the length list of a globally optimal path partition must include the longest
path length.  Hence $\peb{T}$ equals $\peb{T,r}$ for some peripheral vertex $r$.

We show next that the procedure produces the same length list from each
peripheral vertex.  When $r$ and $r'$ are the endpoints of a longest path $R$,
we showed that $r$-optimal and $r'$-optimal path partitions have the same
length list.  When $R'$ is another longest path from $r'$, the algorithm
would again produce an $r'$-optimal path partition.  Since the lexicographic
order is linear, all $r'$-optimal path partitions have the same length list.

Since every longest path in a tree contains the center of the tree, if the path
joining two peripheral vertices is not a longest path, then each is an endpoint
of a longest path to one other peripheral vertex.  Hence one can move from one
peripheral vertex to any other by at most two instances of ``move to the
opposite end of a longest path''.  Therefore, all peripheral vertices have the
same optimal length list.

Because an optimal path partition must contain a longest path and hence
must be an $r$-optimal path partition for some peripheral vertex $r$, we
conclude that $\peb{T} = \peb{T,r}$ for each peripheral vertex $r$.
\end{Proof}

\vspace{-12pt}
\section{Pebbling Number of Cycles}\label{pebcycles}

Proving an upper bound on the pebbling number requires showing that each of a
large number of distributions is solvable.  The following lemma restricts the
distributions that need to be considered.  A {\em thread} in a graph $G$ is a
path whose vertices have degree 2 in $G$.

\begin{lemma}[{\rm Squishing Lemma}]
For a vertex $r$ in a graph $G$, there is a non-$r$-solvable distribution of
$\peb{G,r}-1$ pebbles on $G$ such that on each thread not containing $r$,
all pebbles occur on just one vertex or on two adjacent vertices.
\end{lemma}
\begin{Proof}
Let $P$ be a thread in $G$.  If a distribution has pebbles on only one
vertex of $P$ or on only two adjacent vertices of $P$, then we say that
$P$ is {\it squished}.

Let $D$ be a distribution of $\peb{G,r}-1$ pebbles that is not $r$-solvable.
We transform $D$ into a distribution of the same size such that every thread
not containing $r$ is squished.  A {\it squishing move} removes 1 pebble from
each of two vertices on a thread and puts 2 pebbles on some vertex between them
on the thread.  If some path $P$ is not squished, then we can perform a
squishing move on $P$.  Each squishing move reduces the value of
$\sum_p 2^{-b(p)}$, where the sum is over the set of pebbles on $P$ and $b(p)$
is the distance of pebble $p$ from a fixed end of $P$.  Thus a sequence of
squishing moves must end by squishing $P$.

Let $D'$ be the result of a squishing move applied to $D$ on a thread of $P$ not
containing $r$; pebbles from $y$ and $z$ are moved to $x$ between them.  We
show that if $D'$ is $r$-solvable, then $D$ is $r$-solvable.  Let $\sigma$ be a
pebbling sequence from $D'$ that reaches $r$.  If $\sigma$ never moves pebbles
off $x$, then $\sigma$ also reaches $r$ from $D$.  Hence we may assume that
$\sigma$ includes a move from $x$ to a neighbor $x'$, which we may assume is
toward $y$ along $P$.

By the No-Cycle Lemma, we may assume that $\sigma$ makes no move from $x'$ to
$x$.  The two pebbles used to move from $x$ to $x'$ thus produce no more
benefit than the one pebble that started on $y$ in $D$; under $D$ starts
farther that $x'$ in the only direction it can go.  Also it cannot hurt to have
the extra pebble on $z$.  Thus $D$ also is $r$-solvable.
\end{Proof}

The Squishing Lemma provides a short proof for the pebbling number of $C_n$.

\begin{theorem}[{\rm Pachter et al~\cite{pachter95}}]
The pebbling number of the cycle satisfies $\peb{C_{2k}}=2^k$ and
$\peb{C_{2k+1}}=2\floor{2^{k+1}/3}+1$.
\end{theorem}

\begin{Proof}
{\it Lower Bound.}
Given a root $r$ in $C_{2k}$, a distribution with $2^k-1$ pebbles on the one
vertex at distance $k$ from $r$ is not $r$-solvable.  We show that in
$C_{2k+1}$, a distribution with $\FL{2^{k+1}/3}$ pebbles on each of the two
vertices at distance $k$ from $r$ is not $r$-solvable.  One pile alone cannot
move distance $k$ to reach $r$.  If we combine them first, moving half of one
pile to the other, then the resulting pile has at most
$\FR{2^{k+1}-1}3+\FR12\FR{2^{k+1}-1}3$ pebbles, since $2^{k+1}$ is not
divisible by 3.  The sum is less than $2^k$, so again the pile cannot reach $r$.

{\it Upper Bound.}
A distribution having $2^k$ pebbles on some path of length $k$ ending $r$
is $r$-solvable, since $\peb{P_{k+1}}=2^k$.  This suffices for most cases,
since the Squishing Lemma allows us to restrict attention to distributions
covering only one or two adjacent vertices.  In $C_{2k}$, every two adjacent
vertices lie together in a path of length $k$ ending at $r$.  This also holds
for all cases in $C_{2k+1}$ except when the two adjacent vertices are the two
vertices $s$ and $s'$ at distance $k$ from $r$.

In this case, with all the pebbles on $\{s,s'\}$, we move as many as possible
from the vertex with fewer pebbles to the vertex with more pebbles.  With $m$
pebbles total and $l$ in the smaller pile, the new pile has size at least
$m-l+\FL{l/2}$.  Since $l\le \FL{m/2}$ and $m\ge 2\FL{2^{l+1}/3}+1$, we
obtain a pile of size at least $2^k$ at distance $k$ from $r$, which suffices.
\end{Proof}

\vspace{-12pt}
\section{Optimal Pebbling Number} \label{optimal-pebbling-section}

For optimal pebbling numbers, upper bounds are generally easier than lower
bounds.  For an upper bound, we give a distribution and show that it is
solvable.  For a lower bound, we must show that every distribution up to a
certain size is not solvable.

The Smoothing Lemma plays the role for optimal pebbling that the Squishing Lemma
plays for ordinary pebbling.  The purpose again is to restrict the form of
distributions we study to determine the value of the parameter.  Instead of
squishing pebbles together on a thread, we spread them out.

When $D$ is a distribution on a graph with a vertex $v$ of degree 2, and $v$
has at least three pebbles in $D$, a {\it smoothing move} from $v$ changes $D$
by removing two pebbles from $v$ and adding one pebble at each neighbor of $v$.
The case $m=2$ below will be used in Section~\ref{secmin}.

\begin{lemma}\label{pre-smoothing}
Let $D$ be a distribution on a graph $G$ with distinct vertices $u$ and $v$,
where $v$ has degree 2.  If $D(v)\ge3$, and $u$ is $m$-reachable under $D$,
then $u$ is $m$-reachable under the distribution $D'$ obtained by making a
smoothing move from $v$.
\end{lemma}

\begin{Proof}
For any pebbling sequence $\sigma$ starting from $D$, we form a sequence
$\sigma'$ from $D'$.  If $\sigma$ never makes a move from $v$, then we may set
$\sigma'=\sigma$, since at each step there are at least as many pebbles at each
vertex other than $v$ when starting with $D'$.

If $\sigma$ makes a move from $v$, then let $\sigma'$ be the same as $\sigma$
except that $\sigma'$ skips the first such move.  Having made that move,
$\sigma$ on $D$ produces the same configuration as $\sigma'$ on $D'$, except
that $\sigma'$ on $D'$ has an extra free pebble on one neighbor of $v$.  We
complete $\sigma'$ using the rest of $\sigma$ and have the same number of
pebbles at each vertex as under $\sigma$ from $D$, plus an extra pebble on one
neighbor of $v$.  (Since $\sigma'$ mimics $\sigma$, we never use that extra
pebble.)
\end{Proof}

A distribution $D$ is {\it smooth} if it has at most two pebbles on every
vertex of degree 2 (so no smoothing move is possible).  A vertex $D$ is {\it
unoccupied} under $D$ if $D(v)=0$.

\begin{lemma}[{\rm Smoothing Lemma}] \label{smoothing}
If $G$ is connected and $n(G)\geq 3$, then $G$ has a smooth minimal solvable
distribution with all leaves unoccupied.
\end{lemma}

\begin{Proof}
A minimal solvable distribution has $\optpeb{G}$ pebbles, and always
$\optpeb{G}\le n(G)$.  We first transform an arbitrary solvable distribution
$D$ with $|D|\le n(G)$ into a smooth solvable distribution of the same size;
later we also eliminate pebbles from leaves.

By Lemma~\ref{pre-smoothing}, a smoothing move from $v$ preserves the
reachability of vertices other than $v$.  Since a smoothing move from $v$
leaves a pebble at $v$, also $v$ remains reachable.  Therefore, smoothing
moves preserve solvability.  To complete the proof of the first claim, it
suffices to show that a smooth distribution will result from applying smoothing
moves to any distribution with at most $n(G)$ pebbles.

Suppose first that $G$ is not a cycle.  Starting from any distribution on $G$,
we show that only finitely many smoothing moves can be made.  Every vertex $v$
of degree 2 lies in a unique maximal thread.  Let $P$ be the unique path
through $v$ whose internal vertices have degree 2 and whose endpoints do not.
When $P$ has length $m$ and $v$ has distance $k$ from one end of $P$, we count
each pebble on $v$ with weight $k(m-k)$; it does not matter which end the
distance is measured from.  Pebbles on a vertex with degree other than 2 count
with weight 0.

Let $v$ be a vertex at distance $k$ from the end of a thread of length $m$
(here the ends have degree other than 2).  A smoothing move from $v$ replaces
weight $2k(m-k)$ at $v$ with weight $(k-1)(m-k+1)+(k+1)(m-k-1)$ at its
neighbors.  The total weight declines by $2$.  It must remain nonnegative, so
we reach a distribution with no smoothing move available.

When $G$ is a cycle, we use induction on the number of unoccupied vertices.
Since $|D|\le n(G)$, when all vertices are occupied there is one pebble on
each vertex and $D$ is smooth.  If $v$ has no pebbles and $D$ is not smooth,
then we view $v$ as both endpoints of a thread around the cycle.  Using the
same weight argument as above, each smoothing move reduces the total weight by
2.  Thus eventually the distribution becomes smooth or a pebble moves to $v$.
Since smoothing never uncovers a vertex, moving a pebble to $v$ reduces the
number of unoccupied vertices.  Thus the continuation of the smoothing process
produces a smooth distribution.

We have obtained a smooth minimal solvable distribution $D$; now we consider
leaves.  Let $v$ be a leaf, and let $u$ be its neighbor.  Suppose
that $D(u)=j$ and $D(v)=k\ge1$.

{\em Case 1: $j+k\ge3$}.
Modify $D$ by deleting the pebbles on $v$ and adding $k-1$ pebbles to $u$
instead.  The resulting $D'$ is still solvable, since $D'(u)\ge2$ makes $v$
reachable, and $D'$ starts with at least as many pebbles on $u$ as $v$ could
send there to help pebble other vertices.  However, $|D'| < |D|$, which
contradicts the minimality of $D$.

{\em Case 2: $j+k=2$}.
Modify $D$ by putting both pebbles on $u$.  Still $D'$ is smooth if $u$ has
degree 2.  The two pebbles can be used to cover $v$, and they provide as much
help for other vertices as before.

{\em Case 3: $(j,k)=(0,1)$}.
Move the one pebble to $u$; again $D'$ is smooth.  Because $D$ is $u$-solvable
and cannot use the pebble on $v$ to reach $u$, we can now move another pebble
to $u$ and use the two of them to reach $v$.
\end{Proof}

The Smoothing Lemma yields a short proof of the result of Pachter
et al~\cite{pachter95} that $\optpeb{P_n}=\ceil{2n/3}$, and it yields the
same value also for cycles.  Another short proof was given by
Friedman and Wyels~\cite{friedman}.  We separate an observation
useful in Section~\ref{girth-sec}.

\begin{lemma}\label{pathend}
Let $v$ be an unoccupied vertex in a smooth distribution $D$ on a path with
at most two pebbles on each endpoint.  If $v$ is an endpoint, then $v$ is not
$2$-reachable under $D$.  If $v$ is an internal vertex, then no pebbling
sequence can move a pebble out of $v$ without using an edge in both directions.
\end{lemma}
\begin{Proof}
The first claim follows immediately from the case $m=2$ of the Weight Argument,
since each vertex has at most two pebbles.  For the second claim, moving a
pebbling out of $v$ without first moving a pebble in from each neighbor would
require contradicting the first claim on a smaller path.
\end{Proof}

\vspace{-12pt}
\begin{theorem} \label{cycle-optpeb-theorem}
$\optpeb{C_n}=\optpeb{P_n} = \ceil{2n/3}$.
\end{theorem}
\begin{Proof}
Let $G$ be $C_n$ or $P_n$.

{\it Upper Bound}.
Partition $G$ into $\floor{n/3}$ copies of $P_3$ and possibly one or two
leftover vertices.  Put two pebbles on the central vertex of each $P_3$ and
one pebble on each of the leftover vertices (if any exist).  The distribution
is solvable and has size $\ceil{2n/3}$.

{\it Lower Bound}.
By Lemma \ref{smoothing}, it suffices to consider a smooth solvable
distribution $D$ with no pebbles on leaves.  We use induction on $n$, checking
$n\leq 5$ exhaustively.

By the No-Cycle Lemma, we may assume that the directed edges representing
moves in a pebbling sequence to reach a target vertex form edge-disjoint paths,
and no edge is used in both directions.  Since $D$ is smooth,
Lemma~\ref{pathend} implies that each such path has no unoccupied internal
vertex.


Since $\optpeb{G}\leq\ceil{2n/3}$ and $n\ge6$, at least two vertices of $G$
are unoccupied.  We may choose three unoccupied vertices, since otherwise
$n=6$, no vertex has two pebbles, and $D$ is not solvable.  With three
unoccupied vertices, we can choose an unoccupied internal vertex in $P_n$
or nonadjacent unoccupied vertices in $C_n$; let $S$ be this chosen set.

Since pebbles cannot be sent across an unoccupied vertex, $S$ splits $G$ into
two paths, each of which cannot contribute pebbles to help pebble a vertex on
the other path.  Since the distribution is solvable, each vertex of $S$ can be
pebbled; we treat the vertex as being part of the path that pebbles it,
choosing one such path if both can pebble it.

We now have paths of order $l$ and $n-l$ with $1\le l\le n-1$, and $D$ breaks
into solvable distributions for these two paths.  By the induction hypothesis,
the number of pebbles in $D$ is at least $\ceil{2l/3}+\ceil{2(n-l)/3}$, which
is at least $\ceil{2n/3}$.
\end{Proof}

Next we show that the path is a hardest tree for optimal pebbling number.
It is far from unique; there are many trees whose optimal pebbling number is
$\ceil{2n/3}$.  We write $d(v)$ for the degree of a vertex $v$, and $N(v)$ for
the set of vertices adjacent to $v$.

\begin{theorem}
If $T$ is an $n$-vertex tree, then $\optpeb{T}\le \ceil{2n/3}$.
\label{upperbound}
\end{theorem}
\begin{Proof}
We use induction on $n$.  The claim holds for $n\le 3$, since all such trees
are paths.  In the induction step ($n>3$), we delete three or more vertices at
or near the end of a longest path in $T$ to obtain a subtree $T'$.  It suffices
to show that we can add two pebbles to a minimal solvable distribution $D'$ on
$T'$ to form a solvable distribution $D$ on $T$.  When we add pebbles to $D'$,
all vertices in $T'$ remain reachable, so the problem reduces to showing that
the new vertices can be reached.

Let $P$ be a longest path in $T$.  Let $z$ be an endpoint of $P$, adjacent to
$y$, and let $x$ be the other neighbor of $y$ on $P$.  We consider four cases.

{\em Case 1: $d(y)>2$}.
Since $P$ is a longest path, all neighbors of $y$ other than $x$ are leaves.
Let $T'=T-y-(N(y)-\{x\})$.  Form $D$ from $D'$ by adding two pebbles on $y$;
these make leaf neighbors of $y$ reachable.

{\em Case 2: $d(x)=d(y)=2$}.
Let $T'=T-\{x,y,z\}$.  Form $D$ from $D'$ by adding two pebbles on $y$;
these make $x$ and $z$ reachable.

{\em Case 3: $d(y)=2$ and $x$ has a leaf neighbor $u$}.
Let $T'=T-\{u,y,z\}$.  Form $D$ from $D'$ by adding two pebbles on $y$.
Now $y$ and $z$ are reachable.  We can also reach $u$ by moving a pebble to
$x$ using the distribution $D'$ on $T'$ and then moving a second pebble to $x$
from $y$.

{\em Case 4: $d(y)=2$ and $x$ has no leaf neighbors}.
Let $u$ be a neighbor of $x$ outside $P$.  Since $P$ is a longest path, every
neighbor of $u$ other than $x$ is a leaf.  Let $v$ be a leaf neighbor of $u$,
and let $T'=T-\{v,y,z\}$.
If $x$ is $2$-reachable under $D'$, then we form $D$ by adding two pebbles on
$x$, making $\{v,y,z\}$ all reachable under $D$.  If $u$ is $2$-reachable under
$D'$, then $v$ is reachable, so we form $D$ by adding two pebbles on $y$.  If
neither $x$ nor $u$ is $2$-reachable under $D'$, then no pebbling sequence
starting with $D'$ uses the edge $xu$ in either direction.  Hence from $D'$ we
can reach $x$ and $u$ simultaneously.  Now we form $D$ by adding two pebbles
on $y$, making $\{v,y,z\}$ all reachable after moving pebbles to both $x$ and
$u$ using $D'$.
\end{Proof}

\vspace{-18pt}

\begin{corollary}\label{n-upper}
If $G$ is a connected $n$-vertex graph, then $\optpeb{G}\leq\ceil{2n/3}$,
which is sharp.
\end{corollary}

\vspace{-3pt}
\begin{Proof}
Adding an edge to a graph cannot increase its optimal pebbling number.  Since
$G$ is connected, it has a spanning tree $T$.  Applying
Theorem~\ref{upperbound} to $T$ gives the bound, which is achieved by $P_n$.
\end{Proof}

\vspace{-9pt}
Finally, we give a short proof that $\optpeb{Q_k} \geq (4/3)^k$.  The proof by
Moews~\cite{moews98} used a continuous relaxation of pebbling, but the standard weight function and expectation suffice.

\begin{theorem}[{\rm Moews~\cite{moews98}}]
$\optpeb{Q_k}\geq\left(4/3\right)^{k}$, where $Q_k$ is the $k$-dimensional
hypercube.
\end{theorem}

\vspace{-3pt}
\begin{Proof}
Let $D$ be a solvable distribution on $Q_k$; we show that $|D|\geq (4/3)^k$.
Since $D$ is solvable, the standard weight inequality
$\sum_{i\ge0} a_{i,r}2^{-i}\ge1$ holds for each vertex $r$, where
$a_{i,r}$ is the number of pebbles at distance $i$ from $r$ in $D$.

Select a vertex $r$ in $Q_k$ uniformly at random.  Since the weight inequality
holds for each $r$, linearity of expectation yields
$\sum_{i\geq0}2^{-i}\mathbf{E}\left[a_{i,r}\right]\geq 1$.
For a fixed pebble on a vertex $u$, the probability that $r$ has
distance $i$ from $u$ is ${k\choose i}2^{-k}$, since $Q_k$ has
$2^k$ vertices and ${k \choose i}$ of them have distance $i$ from $r$.
By linearity of expectation,
$ \mathbf{E}\left[a_{i,r}\right]=|D|{k\choose i}2^{-k}$.
Substituting and simplifying now yields

\vspace{-9pt}

\begin{eqnarray*}
|D|\sum_{i\geq0}{k\choose i}2^{-i} & \geq & 2^{k}.
\end{eqnarray*}
Applying the Binomial Theorem yields $|D|(1+\frac12)^k\ge 2^k$, and hence
$|D|\ge(4/3)^k$.
\end{Proof}

\vspace{-12pt}
\section{Bounds in Terms of Minimum Degree} \label{secmin}

We have proved that $\optpeb{G} \leq \ceil{2n/3}$ for every connected
$n$-vertex graph $G$, with equality for paths and cycles.  One would expect
that tighter upper bounds hold for denser graphs.  How large can $\optpeb{G}$
be when we require minimum degree $k$?

A {\em dominating set} in a graph $G$ is a set $S\esub V(G)$ such that every
vertex not in $S$ has a neighbor in $S$.  The {\em domination number}
$\gamma(G)$ is the minimum size of a dominating set.  Placing two pebbles at
each vertex of a dominating set yields $\optpeb{G}\le 2\gamma(G)$.  Thus
upper bounds on $\gamma(G)$ yield upper bounds on $\optpeb{G}$.

For graphs with minimum degree at least $k$, Arnautov~\cite{arnau} and
Payan~\cite{payan} proved that $\gamma(G)\le n\frac{1+\ln(k+1)}{k+1}$; a short
probabilistic argument appears in Alon~\cite{alon}.  In a $k$-regular
$n$-vertex graph, dominating sets have size at least $\frac n{k+1}$, and
Alon~\cite{alon} showed that the domination number may be as large as
$(1+o(1))n\FR{1+\ln(k+1)}k$.  Hence we cannot improve the bound using
domination number alone.

Czygrinow~\cite{czygrinow} communicated to us an easy argument for a better
upper bound when $k\ge3$; we begin by presenting this.  A {\em distance-2
dominating set} in a graph $G$ is a set $S\esub V(G)$ such that every vertex of
$G$ is within distance at most 2 from $S$ (similarly, one can define
distance-$d$ dominating sets).  The case $d=1$ of the following proposition is
folklore in some circles but seems to be unknown in the subject of graph
domination.  We will use the general result in Section~\ref{girth-sec}.

\begin{proposition}\label{distd}
If $c$ is the minimum size of a distance-$d$ neighborhood in $G$, then
$G$ has a distance-$2d$ dominating set of size at most $n(G)/c$.
\end{proposition}
\begin{Proof}
We build such a set $S$.  Initially, put one vertex in $S$.  As we proceed, let
$T$ consist of all vertices within distance $d$ of $S$.  If $T$ is not a
distance-$d$ dominating set, then let $v$ be a vertex that is not within
distance $d$ of $T$.  Add $v$ to $S$; this adds the distance-$d$ neighborhood
of $v$ to $T$, none of which was in $T$ before.  Thus $T$ grows by at least $c$
vertices for each vertex added to $S$.  We therefore add at most $n/c$ vertices
to $S$ by the time $T$ becomes a distance-$d$ dominating set, at which point
$S$ is a distance-$2d$ dominating set.
\end{Proof}

\vspace{-12pt}

\begin{corollary}[{\rm Czygrinow}] \label{mindeg-upper}
If $G$ is a graph with minimum degree $k$, then
$\optpeb{G}\le \frac{4n(G)}{k+1}$.
\end{corollary}
\begin{Proof}
Distance-1 neighorhoods have size at least $k+1$, so
Proposition~\ref{distd} yields a distance-2 dominating set $S$ of
size at most $n(G)/(k+1)$.  Put four pebbles at each vertex of $S$.
\end{Proof}

Corollary~\ref{mindeg-upper} improves the upper bound of $\CL{2n/3}$
from Corollary~\ref{n-upper} when $k\ge6$.  A simple construction shows
that this easy upper bound is within a factor of 2 of being sharp; we present
$n$-vertex graphs with minimum degree $k$ whose solvable distributions have at
least $\frac{2n}{k+1}$ pebbles.  Subsequently we present a better construction
with optimal pebbling number approximately $\frac{2.4n}{k+1}$.

We begin by introducing another technique for proving lower bounds.
Given a graph $G$, the operation of {\em collapsing} a vertex set
$S\subseteq V(G)$ produces a new graph $H$ in which $S$ is replaced with
a single vertex whose neighbors are the neighbors of $S$ in $G$ that were
outside $S$.  The subgraph induced by $V(G)-S$ remains unchanged.
We use the term ``collapsing'' rather than ``contracting'' because the
subgraph of $G$ induced by $S$ need not be connected.

\begin{lemma}[{\rm Collapsing Lemma}]\label{collapsing}
If $H$ is obtained from $G$ by collapsing vertex sets,
then $\optpeb{G}\ge\optpeb{H}$.
\end{lemma}
\begin{Proof}
Let $D$ be a solvable distribution on $G$.  Form distribution $D'$ on
$H$ as follows: for each collapsed set $S$, put all the pebbles that
were on $S$ in $D$ onto the single vertex representing $S$ in $H$.
Treat uncollapsed vertices as collapsed sets of size 1.

To show that $D'$ is solvable, for $u\in V(H)$ choose a vertex $v\in V(G)$ in
the set that collapses to $u$.  Let $\sigma$ be a pebbling sequence
from $D$ that reaches $v$.  The sequence $\sigma$ ``collapses'' in an
obvious way to a sequence $\sigma'$ from $D'$ that reaches $u$.
More precisely, the distribution $C$ resulting from a pebbling move on $D$
collapses to a distribution $C'$ on $H$ that is obtained from $D'$ by
discarding one pebble (if the move on $D$ was within a collapsed set) or by
making one pebbling move from $D'$.
\end{Proof}

\vspace{-12pt}

\begin{proposition}\label{clique-ring}
For $n>k\ge2$, there is an $n$-vertex graph $G$ with minimum degree $k$ such
that $\optpeb{G}>\FR{2n}{k+1}-2$, improving to $\optpeb{G}\ge\FR{2n}{k+1}$
when $n$ is a multiple of $k+1$.
\end{proposition}
\begin{Proof}
When $n=k+1$, the complete graph $K_n$ has this behavior.

When $n$ is a larger multiple of $k+1$, let $J$ be the graph obtained from
$K_{k+1}$ by deleting one edge; the {\it internal} vertices of $J$ are its
vertices of degree $k$.  Let $G$ be the $k$-regular ``ring of cliques'' with
$r(k+1)$ vertices formed by putting $r$ copies of $J$ in a circle and making
one non-internal vertex in each copy of $J$ adjacent to one non-internal
vertex in the next copy.

By Lemma~\ref{collapsing}, collapsing the internal vertices in a copy of
$J$ into one vertex cannot increase the optimal pebbling number.  Doing this
in each copy of $J$ produces $C_{3r}$.  By Theorem~\ref{cycle-optpeb-theorem},
we obtain $\optpeb{G} \geq 2r=2n/(k+1)$.

For general $n$, let $r=\FL{n/(k+1)}$.  Form $J'$ from $K_{n-(r-1)(k+1)}$ by
deleting one edge.  Form $G'$ by the construction for $G$ above, using one copy
of $J'$ and $r-1$ copies of $J$.  Collapsing $n+1-r(k+1)$ internal vertices of
$J'$ into one vertex turns $G'$ into the example $G$ for $r(k+1)$ vertices.
By Lemma~\ref{collapsing}, $\optpeb{G'}\ge\optpeb{G}\ge 2r\ge 2(n-k)/(k+1)$.
\end{Proof}

\vspace{-9pt}
Corollary~\ref{n-upper} shows that the construction in
Proposition~\ref{clique-ring} is extremal for $k=2$, where it produces $C_n$.
For $k=3$, it provides connected $n$-vertex graphs with optimal pebbling number
asymptotic to $n/2$; the upper bound from Corollary~\ref{n-upper} remains
$\CL{2n/3}$.  As $k$ grows, the coefficient on $n$ in
Proposition~\ref{clique-ring} decreases.

However, for $k>15$ the optimal pebbling number of our next construction
exceeds $2\frac n{k+1}$ asymptotically for large $n$.  In particular, there is
an $n$-vertex graph $G_n$ with minimum degree $k$ such that
$\optpeb{G_n}\FR{k+1}n\to 2.4-\FR{24}{5k+15}$.  This limit exceeds 2 when
$k>15$.  We present the construction only for $k\equiv 0\mod{3}$; slightly
weaker results hold for general $k$.

We will apply Lemma~\ref{collapsing} to a graph that we will contract to a
cycle.  We first develop a lower bound for 2-solvable distributions on cycles.

\begin{lemma} \label{smooth2}
Let $G$ be a graph with distribution $D$, and let $A$ be a subset of $V(G)$
such that each vertex in $A$ has a neighbor in $A$.  If each vertex in $A$ is
$2$-reachable under $D$, then each vertex in $A$ is $2$-reachable under any
distribution produced from $D$ by a smoothing move.
\end{lemma}
\begin{Proof}
Let $D'$ be a distribution obtained from $D$ by a smoothing move from $v$.
Note that $D'(v)\ge1$, by the definition of smoothing.  By
Lemma~\ref{pre-smoothing}, every vertex of $A-\{v\}$ is $2$-reachable
under $D'$.  Hence we may assume that $v\in A$.

Let $u$ be a neighbor of $v$ in $A$, and let $\sigma$ be a pebbling sequence
under $D'$ after which $u$ has two pebbles.  If $\sigma$ has a move out of $v$,
then truncating $\sigma$ yields a pebbling sequence showing that $v$ is
$2$-reachable.  Otherwise, $v$ retains at least one pebble after executing
$\sigma$, and then a pebbling move from $u$ to $v$ gives it another.
\end{Proof}

\vspace{-18pt}
\begin{lemma} \label{cycl2solv}
For $n\ge3$, if at least $n-1$ vertices are $2$-reachable under a distribution
$D$ on $C_n$, then $|D|\geq n$.
\end{lemma}
\begin{Proof}
Having a $2$-reachable vertex requires that $D$ has two pebbles on some vertex.
This completes the proof when there is at most one unoccupied vertex.  Hence we
may choose distinct unoccupied vertices $u$ and $v$.  With $n-1$ vertices
$2$-reachable under $D$, Lemma~\ref{smooth2} and the weight argument used in the
Smoothing Lemma allow us to assume that $D$ is smooth.

\nobreak
Let $P$ and $P'$ be the $u,v$-paths along the cycle.  Since at least $n-1$
vertices are $2$-reachable, we may assume that $u$ is $2$-reachable.  By
Lemma~\ref{pathend}, a pebbling sequence cannot move a pebble out of $v$ without
using an edge in both directions, which by the No-Cycle Lemma
does not occur in some pebbling sequence that moves two pebbles to $u$.
Lemma~\ref{pathend} also implies that $u$ is not $2$-reachable under the
restrictions of $D$ to $P$ or $P'$.  Therefore, $2$-reachability of $u$
requires moving a pebble to $u$ from each of $P$ and $P'$, independently.
Hence each must have a vertex with two pebbles.

In particular, there is a vertex with two pebbles on each path of occupied
vertices joining two unoccupied vertices, and therefore $|D|\ge n$.
\end{Proof}

Let $G_{r,s}$ be a graph formed from $r$ disjoint copies of $K_s$ in a row
by making each vertex adjacent to all but one vertex in each neighboring
copy of $K_s$ (there is only one such graph, up to isomorphism).  Similarly,
let $H_{r,s}$ be a graph formed from $r$ disjoint copies of $K_s$ in a circle
via the same definition (the isomorphism class is determined by the placement
of edges joining the last two copies of $K_s$).

\begin{theorem}\label{4r5}
For $s\ge 3$ and $r\ge 1$, $\CL{4r/5}\le \optpeb{G_{r,s}} \le 4\CL{r/5}$.
If $r,s\ge 3$, then $\CL{4r/5}\le \optpeb{H_{r,s}} \le 4\CL{r/5}$.
The lower bounds hold also when $s=2$.
\end{theorem}
\begin{proof}
Call the initial copies of $K_s$ the ``cliques''.  By dividing the cliques into
consecutive groups of five and placing four pebbles on some vertex in the
central clique of the group, we obtain a solvable distribution that uses
$4\CL{r/5}$ pebbles.  In particular, note that if $s\ge3$, then any two
vertices in cliques that are two apart in the ring have a common neighbor in
the intervening clique.  This fails when $s=2$, and hence the upper bounds
require $s\ge3$ (the lower bound for $s=2$ is strengthened in
Theorem~\ref{ladder}).

Now consider the lower bounds.
The proof is by induction on $r$.  When $r\le 4$, the claims are easily checked.
Since adding edges to a graph cannot increase the optimal pebbling number,
it suffices in the induction step to prove the lower bound for $H_{r,s}$.

For $r \geq 5$, consider a minimal solvable distribution $D$ on $H_{r,s}$.
Label the cliques $F_1, F_2, \ldots, F_r$ in order.  Let $A$ be the set of
cliques containing no vertex that is $2$-reachable under $D$.  If $|A| \leq 1$,
then collapsing each $F_i$ to a single vertex yields a distribution on $C_r$
under which at least $r-1$ vertices are $2$-reachable.
By Lemma~\ref{cycl2solv}, in this case $|D| \ge r > \CL{4r/5}$.

We may therefore assume that $|A|\ge2$.  Suppose first that $F_i \in A$ but
$F_{i-1}, F_{i+1} \notin A$.  Let $u$ be a $2$-reachable vertex in $F_{i-1}$,
and let $v$ be a $2$-reachable vertex in $F_{i+1}$.  Since $|A| \geq 2$, we may
also choose $F_j \in A$.  Since $F_i,F_j\in A$, we can never put two pebbles
on any vertex in $F_i\cup F_j$, and hence we can never move a pebble out of
$F_i\cup F_j$.  Since $u$ and $v$ are separated by $F_i\cup F_j$, this implies
that $u$ and $v$ are $2$-reachable simultaneously; that is, the pebbles used in
moving two pebbles to one are not used in moving two pebbles to the other.
Since $s \geq 3$, $u$ and $v$ have a common neighbor $w$ in $F_i$.  Now $w$ is
$2$-reachable using pebbles moved from $u$ and $v$, which contradicts
$F_i\in A$.

It follows that for every member of $A$, some neighboring clique is also in $A$.
When $F_i,F_{i+1}\in A$, we call the edges joining $F_i$ and $F_{i+1}$
{\it useless}.  Since we cannot move two pebbles to any vertex in either
clique, we cannot move a pebble along an edge joining them.  Hence deleting
these edges does not affect the solvability of $D$.

Since for every member of $A$ there is a neighboring clique also in $A$, every
clique in $A$ is incident to a useless set of edges.  Hence there are at least
$|A|/2$ such useless sets of edges.

If $|A|\ge 3$, then there are at least two useless sets of edges; deleting
them leaves a graph whose components are $G_{t,s}$ and $G_{r-t,s}$, with
the distribution $D$ still solvable.  Applying the induction hypothesis to the
two components yields $|D|\geq(4/5)r$.

Otherwise, $|A| = 2$.  Lemma~\ref{collapsing} implies that collapsing each
clique to a single vertex and collapsing the two vertices arising from $A$ to a
single vertex $v$ yields a distribution on $C_{r-1}$ under which every vertex
except $v$ is $2$-reachable.  Since its size is $|D|$, Lemma~\ref{cycl2solv}
implies that $|D|\ge r-1\ge 4r/5$.
\end{proof}

\vspace{-15pt}
\begin{corollary}\label{mindeg-lower}
Let $k$ be a positive multiple of 3.
For $n\ge k+3$, there is an $n$-vertex graph $G$ with minimum degree $k$
such that $\optpeb{G}\ge(2.4-\frac{24}{5k+15}-o(\FR1n))\frac n{k+1}$.
When $n$ is a multiple of $(k/3)+1$, the term $-o(\FR1n)$ can be dropped.
\end{corollary}

\vspace{-3pt}
\begin{Proof}
Given such $n$ and $k$, let $s=k/3+1$ and $r=\FL{n/s}$.  Note that $r,s\ge 3$.
The graph $H_{r,s}$ is $3(s-1)$-regular, since each vertex has $s-1$ neighbors
in its own clique and in each neighboring clique.  Form $G$ by adding to
$H_{r,s}$ a set of $n-rs$ vertices whose neighborhoods duplicate neighborhoods
of vertices in $H_{r,s}$.  Thus $G$ has $n$ vertices and minimum degree at
least $k$.  Also $H_{r,s}$ is obtained from $G$ by collapsing sets of vertices.
Thus $\optpeb{G}\ge\optpeb{H_{r,s}}$, by Lemma~\ref{collapsing}.
When $n$ is a multiple of $s$, we compute
$$
\frac{\optpeb{G}(k+1)}{n}=
\frac{\optpeb{H_{r,s}}(3s-2)}{rs}\ge
\frac{4r}5\frac{3s-2}{rs}=
\frac{12}5-\frac{8}{5s}=
\frac{12}5-\frac{24}{5k+15}.
$$
In general, $n\le rs+s-1$, so we replace the $rs$ in the denominator above
with $rs(1+\FR{s-1}{rs})$.  Since
$(1+\FR{s-1}{rs})^{-1}\ge 1-\FR{s-1}{rs}\ge 1-\FR{k}{3n-k}$, we obtain
$\FR{12}5-\FR{24}{5k+15}-o(\FR1n)$ as a lower bound.
\end{Proof}

Let $f(k)$ be the infimum of all $\alpha$ such that
$\optpeb{G}\FR{k+1}{n(G)}\le \alpha$ for all graphs with minimum degree $k$.
By Corollary~\ref{mindeg-upper}, Proposition~\ref{clique-ring}, and
Corollary~\ref{mindeg-lower}, $\max\{2,2.4-\frac{24}{15k+5}\}\le f(k)\le 4$.
Given the simplicity of Corollary~\ref{mindeg-upper}, we believe that
$f(k)$ is bounded away from 4, but we have no conjecture for an asymptotic
value.

For $k=3$, the upper bound $\optpeb{G}\le 2n(G)/3$ yields $f(3)\le 8/3$.
We have no construction needing more than the $n(G)/2$ of
Proposition~\ref{clique-ring}; Theorem~\ref{ladder} provides another
such example.

\begin{question}
Is it true that $\optpeb{G}\le \CL{n/2}$ whenever
$G$ is a connected $n$-vertex graph with minimum degree at least 3?
The bound would be sharp for $n\ge6$.
\end{question}

\bigskip
When $k=4$, Corollary~\ref{mindeg-lower} does not apply, but more pebbles may
be neeeded than the $2n/5$ in Proposition~\ref{clique-ring}.  We base our
construction on the ``Sierpinski Triangle''.

\begin{example}
{\rm
Let $G_1$ be a triangle; its three vertices are its {\it corners}
$\{x,y,z\}$.  For $m>1$, given three copies of $G_{m-1}$ with corner vertices
$\{x_i,y_i,z_i\}$ in the $i$th copy, form $G_m$ by collapsing the pairs
$\{z_1,x_2\}$, $\{y_2,z_3\}$, and $\{x_3,y_1\}$.  The remaining corner vertices
$\{x_1,y_2,z_3\}$ are the corners of $G_m$.  Another way to construct
$G_m$ from $G_{m-1}$, starting with a layout of $G_1$ in the plane,
is to subdivide the edges of each bounded triangle and add a new triangle
joining each such set of three new vertices.

For $m>1$, form $H_m$ from $G_m$ by adding three edges to make the corners
pairwise adjacent.  Since the corners of $G_m$ have degree 2 and all other
vertices of $G_m$ have degree 4, $H_m$ is 4-regular for $m>1$.  Also,
$n(H_m)=n(G_m)=3n(G_{m-1})-3$; with $n(G_1)=3$, we have $n(H_m)=(3^m+3)/2$.

For $m\ge3$, we present a solvable distribution on $H_m$ with $2\cdot 3^{m-2}$
pebbles (there are many such distributions), and we conjecture that this is
optimal.  If so, then $\optpeb{H_m}/n(H_m)$ approaches $4/9$ from below.

In forming $G_m$, three copies of $G_{m-1}$ are used.  Further breakdown
shows that $3^{m-3}$ copies of $G_3$ are used.  The number $a_m$ of vertices of
$G_m$ that are corners of copies of $G_3$ equals $n(G_{m-2})$, by the
alternative construction.  Since the corners of $G_3$ form a distance-2
dominating set of $G_3$, we have
$\optpeb{H_m}\le \optpeb{G_m}\le 4n(G_{m-2})=2\cdot 3^{m-2}+6$.

For $m\ge3$, we can save six pebbles in this solvable distribution on $H_m$.
The distance between corners of $G_m$ is $2^{m-1}$.  In $H_m$, these corners
are pairwise adjacent.  Hence the four pebbles on one corner $x$ can satisfy
the other corners $y$ and $z$ and the immediate neighbors of $y$ and $z$.
Let $P$ be the shortest $y,z$-path.  If we delete the pebbles on $P$, then the
unreachable vertices are within distance 1 of $P$.  By putting two pebbles
each on the corners of copies of $G_2$ along $P$, we have deleted
$4(2^{m-3}+1)$ pebbles and added $2(2^{m-2}-1)$ pebbles, saving 6.
}
\qed
\end{example}

\section{Girth and Minimum Degree}\label{girth-sec}

Forbidding short cycles restricts the input in a way that improves upper bounds
on the optimal pebbling number.  In particular, if $G$ has minimum degree $k$
and girth at least 5, then four pebbles at a vertex $v$ can take care of
$k^2+1$ vertices, because the neighborhoods of the neighbors of $v$ overlap
only at $v$.

\begin{proposition}\label{girthprop}
If $G$ is a connected graph with minimum degree $k$ and girth at least $2t+1$,
then $\optpeb{G}\leq 2^{2t}n/c_k(t)$, where $c_k(t)=1+k\SE i1t(k-1)^{i-1}$.
\end{proposition}
\begin{Proof}
When $G$ has minimum degree $k$ and girth at least $2t+1$, every distance-$t$
neighborhood has size at least $c_k(t)$.  Proposition~\ref{distd} then applies.
\end{Proof}

Note that $c_k(t)=1+[(k-1)^t-1](1+\frac2{k-2})>(k-1)^t$ for fixed $k$.
For fixed $k$ with $k\ge 6$, this yields $\optpeb{G}/n(G)\to0$
as $t\to\infty$.  A more detailed analysis improves the upper bound.
The idea is to use $2^{2t}$ pebbles on a vertex of the distance-$2t$
dominating set only when it is used to reach substantially more than the
$c_k(t)$ vertices guaranteed in its distance-$t$ neighborhood.

\begin{theorem}\label{biggirth}
Let $k$ and $t$ be positive integers with $k\ge3$ and $t\ge2$, except not
$(k,t)=(3,2)$.
If $G$ is an $n$-vertex graph with minimum degree $k$ and girth at least $2t+1$,
then $\optpeb{G}\leq 2^{2t}n/(c_k(t)+c'(t))$, where $c_k(t)$ is defined as
above and $c'(t)=(2^{2t}-2^{t+1})\frac t{t-1}$.
\end{theorem}
\begin{Proof}
As constructed in the proof of Proposition~\ref{distd}, we begin with a
distance-$2t$ dominating set $S$ of size at most $n/c_k(t)$, where $c_k(t)$
is defined as in Proposition~\ref{girthprop} and the distance between any two
vertices of $S$ is at least $2t+1$.

To each $v\in S$, we assign a set $R(v)$ of vertices in $G$; pebbles on $v$
will be used to reach the vertices of $R(v)$.  Each vertex within distance $t$
of $v$ is in $R(v)$; this causes no conflict, since the distance-$t$
neighborhoods from vertices of $S$ are disjoint.  Indeed, we grow the sets
of the form $R(v)$ to absorb all vertices of $G$ by doing a simultaneous
breadth-first search from all of $S$; each vertex goes into just one of
these sets when it is reached.  Since $S$ is a distance-$2t$ dominating
set, for each $v\in S$ this generates a spanning tree $T(v)$ of the subgraph
induced by $R(v)$, such that leaves of $T(v)$ have distance at most $2t$ from
$v$ in $T(v)$.

Let $R'(v)$ be the set of nonleaf vertices of $T(v)$ that are not within
distance $t$ of $v$.  Let $r'(v)=|R'(v)|$.  If $r'(v)< 2^{2t}-2^{t+1}$, then
put $2^{t+1}$ pebbles on $v$ and one pebble on each vertex of $R'(v)$.
Otherwise, put $2^{2t}$ pebbles on $v$.

When $r'(v)\ge 2^{2t}-2^{t+1}$, the $2^{2t}$ vertices on $v$ can reach
all vertices at distance at most $2t$ from $v$.
When $r'(v)<2^{2t}-2^{t+1}$, the $2^{t+1}$ pebbles on $v$ can reach
vertices at distance $t+1$ from $v$, including the closest ones in $R'(v)$.
The rest of $T(v)$ can then be reached by pebbling along paths through $R'(v)$.
Hence the distribution is solvable.

When $r'(v)<2^{2t}-2^{t+1}$, we use $r'(v)$ pebbles on $R'(v)$.  We claim
that at least $r'(v)\FR t{t-1}$ vertices lie in $T(v)$ that are not within
distance $t$ of $v$.  For $0\le i\le t-1$, let $p_i$ be the number of vertices
in $T(v)$ that are $i$ levels above a leaf, but not within distance $t$ of
$v$.  For $i>0$, the vertices counted by $p_i$ have distinct children in $T(v)$
counted by $p_{i-1}$, so $p_0\ge p_1\ge \cdots \ge p_{t-1}$.  Also,
$r'(v)=\sum_{i=1}^{t-1} p_i$.  We put pebbles on $r'(v)$ vertices, but we add
$r'(v)+p_0$ vertices beyond those counted by $c_k(t)$.  We have
$\FR{r'(v)+p_0}{r'(v)}=1+\FR{p_0}{r'(v)}\ge 1+\FR{p_0}{(t-1)p_0} = \FR t{t-1}$.
Hence we add at least $r'(v)t/(t-1)$ vertices not previously counted.

We have shown that when $r'(v)<2^{2t}-2^{t+1}$, we use $2^{t+1}+r'(v)$ pebbles
with $T(v)$ having at least $c_k(t)+r'(v)\frac t{t-1}$ vertices.
When $r'(v)\ge 2^{2t}-2^{t+1}$, we use $2^{2t}$ pebbles, with
$T(v)$ having at least $c_k(t)+c'(t)$ vertices.

Let $S'=\{v\in S\st r'(v)< 2^{2t}-2^{t+1}\}$, and let $s=|S|$.  Let
$r=\sum_{v\in S'}(2^{2t}-2^{t+1}-r'(v))$.  We have
$n\ge s[c_k(t)+c'(t)]-r\frac t{t-1}$, and we used $2^{2t}s-r$ pebbles.  Thus
$$
\optpeb{G}\le \frac{2^{2t}s-r}{s[c_k(t)+c'(t)]-r\frac t{t-1}}n\le
\frac{2^{2t}}{c_k(t)+c'(t)}n,
$$
where the last inequality uses that $2^{2t}/(c_k(t)+c'(t))<(t-1)/t$ when
$k\ge3$ and $t\ge2$ and $(k,t)\ne(3,2)$.
\end{Proof}

Since $c'(t)\ge 2^{2t}$ and $c_4(t)=1+(4^t-1)(5/3)$, the resulting upper bound
on $\optpeb{G}/n(G)$ when $k=5$ tends to $3/8$ as $t\to\infty$.  For $k=2$,
always $\optpeb{C_n}=\CL{2n/3}$.  Thus it is natural to ask whether the
behavior we noted for $k\ge6$ also holds for $3\le k\le 5$.

\begin{question}
For $k\in\{3,4,5\}$, does there exist $f_k(t)$ such that
$\lim_{t\to\infty}f_k(t)=0$ and graphs with minimum degree $k$ and girth at
least $2t+1$ satisfy $\optpeb{G}/n(G)\le f_k(t)$?
\end{question}

We have not constructed graphs to show that the bound in Theorem~\ref{biggirth}
is sharp, and we do not believe that it is sharp.  We present one more result,
showing that if $G$ has girth 4 and minimum degree 4, then $\optpeb{G}$ can
be as large as $n(G)/2$.  This improves the construction in
Proposition~\ref{clique-ring} for $k=3$ by showing that even when triangles are
forbidden the same number of pebbles may be needed.

The {\it cartesian product} $G\cart H$ of graphs $G$ and $H$ is the graph with
vertex set $V(G)\times V(H)$ such that $(u,v)$ is adjacent to $(u',v')$
if and only if (1) $u=u'$ and $vv'\in V(H)$ or (2) $v=v'$ and $uu'\in E(G)$.
Note that $G\cart H$ contains a copy of $H$ for each vertex of $G$ and a copy
of $G$ for each vertex of $H$.

In particular, $C_m\cart K_2$ and $P_m\cart K_2$ are circular and linear
``ladders''; two copies of the cycle or path, with corresponding vertices
from the two copies adjacent.  We call the $m$ copies of $K_2$ the {\it rungs}
of the graph.  In $C_m\cart K_2$, exchanging the matching joining two rungs
for the other possible matching joining them yields a graph isomorphic to the
graph formed from a $2m$-cycle by adding chords joining opposite vertices
(those at distance $m$ along the cycle).  This graph has been called the
``M\"obius ladder'', so we denote it by $M_m$.

The graphs $C_m\cart K_2$ and $M_m$ are special cases of the construction in
Theorem~\ref{4r5} with $m=r$ and $s=2$.  The lower bound there is $4m/5$; this
result improves that bound.  To prove the lower bound, we need to characterize
the optimal 2-solvable distributions on paths.  For this we need an analogue
of Lemma~\ref{smoothing} for 2-solvable distributions.

\begin{lemma}\label{smoothing2}
Every connected graph with at least three vertices (other than a cycle) has a
smooth minimal 2-solvable distribution that gives at most two pebbles to each
leaf.
\end{lemma}
\begin{Proof}
We apply smoothing moves to a minimal 2-solvable distribution on such a graph
$G$.  Since every vertex is 2-reachable, every vertex has a 2-reachable
neighbor, and hence the result of a smoothing move is also a 2-solvable
distribution, by Lemma~\ref{smooth2}.  We showed in the proof of
Lemma~\ref{smoothing} that when $G$ is not a cycle only finitely many
smoothing moves can be made, so we obtain a smooth minimal 2-solvable
distribution $D$.

Suppose now that $D(v)>2$ for some leaf $v$.  Let $u$ be its neighbor, and
let $j=D(u)$ and $k=D(v)\ge3$.  Obtain $D'$ from $D$ by setting $D'(v)=1$ and
$D'(u)=j+k-2$; leave other values unchanged.  Now $D'$ starts with at least as
many pebbles on $u$ as $v$ could send there under $D$ to help pebble other
vertices.  If $j+k\ge 4$, then $D'(u)\ge2$ to provide a second pebble for $v$.
Otherwise, $(j,k)=(0,3)$; now $v$ can send only one pebble to $u$ under $D$, so
the $2$-solvability of $D$ requires that another pebble can be moved to join
the pebble on $u$ under $D'$; they can then provide a second pebble for $v$.
Hence $D'$ is 2-solvable, but $|D'| < |D|$, which contradicts the minimality of
$D$.
\end{Proof}
%
%
%

A slightly longer case analysis ensures a smooth 2-solvable distribution
with at most one pebble on each leaf, but we will not need this.

\begin{theorem}\label{path2}
Every 2-solvable distribution on $P_n$ has at least $n+1$ pebbles.
Furthermore, the 2-solvable distributions with $n+1$ pebbles consist of
``prime segments'' separated by single unoccupied vertices, where a {\em prime
segment} is a path with either (1) two pebbles on one vertex and one pebble
on all other vertices, or (2) three consecutive vertices having
$0, 4, 0$ pebbles, respectively, and one pebble on all other vertices.
\end{theorem}
\begin{Proof}
We use induction on $n$; when $n\le2$ the unique minimal 2-solvable
distributions have $n+1$ pebbles and are prime segments, as claimed.
Consider $n\ge3$.

By Lemma~\ref{smoothing2}, there is a smooth 2-solvable distribution $D$ having
at most two pebbles on each endpoint.  By Lemma~\ref{pathend}, the endpoints
cannot be unoccupied.  If every vertex is occupied, then 2-solvability requires
some vertex to have two pebbles, and then the minimal distributions have $n+1$
pebbles and form a single prime segment.

We may therefore assume that some internal vertex $v$ is unoccupied.
By Lemma~\ref{pathend}, 2-reachability of $v$ requires one pebble to arrive
from each side.  Since two pebbles cannot arrive at $v$ from one side, pebbles
on one side of $v$ cannot be used to obtain 2-solvability of any
vertex on the other side.  Hence $P_n-v$ consists of two subpaths, each
inheriting a 2-solvable distribution (each neighbor of $v$ is 2-reachable
using only pebbles on that side, because each can provide a pebble to $v$).
With these paths having $l$ and $n-1-l$ vertices, the induction hypothesis
requires $l+1+n-l$ pebbles in $D$, and it also completes the decomposition
into prime segments after the split at $v$.

We now consider other optimal 2-solvable distributions on $P_n$, not
necessarily smooth.  The transformation in Lemma~\ref{smoothing2} shows
that optimal 2-solvable distributions have at most two pebbles on each leaf,
smooth or not.  Since smoothing moves preserve 2-solvability but do not
discard pebbles, a smoothing move on an optimal 2-solvable distribution will
not leave a leaf with at least three pebbles.  Hence we can obtain all
optimal 2-solvable distributions by ``inverting'' smoothing moves starting
with the distributions we have described.

Such an inversion move changes consecutive pebble values $(i,j,k)$ to
$(i-1,j+2,k-1)$, where $i,j,k\ge1$.  Since all values are positive,
we can never make an unoccupied vertex occupied by such a move, so
the three positions must be within a single original prime segment.
We claim that the inversion move maintains the property that 2-solvability
within the segment requires pebbles to flow out from the unique vertex
with most pebbles on the segment, and pebbles never cross an unoccupied
internal vertex.  Maintaining these properties, we can never
make an inversion move with $j=1$, because by symmetry we may assume $i=1$,
and the newly unoccupied vertex would not be 2-reachable.  Hence the only
possible inversion moves change $(1,2,1)$ to $(0,4,0)$, and there can only
be one of these within a prime segment.  Segments formed by surrounding
$(0,4,0)$ with single-pebble vertices are 2-solvable, so this completes the
description of the optimal 2-solvable distributions.
\end{Proof}

\begin{theorem}\label{ladder}
$\optpeb{C_m\cart K_2}= \optpeb{P_m\cart K_2}=\optpeb{M_m}\ge m$ for $m\ge2$.
Equality holds except for $m\in\{2,5\}$.
\end{theorem}
\begin{Proof}
We first provide constructions (except when $m\in\{2,5\}$) to show that the
lower bound is sharp.  Observe that three pebbles on one rung can reach all
vertices on the two neighboring rungs.  Also, four pebbles on two adjacent
rungs (two each at opposite corners of the resulting 4-cycle) can reach all
vertices on the two neighboring rungs.  We can cover the graph with disjoint
sets of three or four rungs unless $m\in\{2,5\}$.  For $m=5$, six pebbles
suffice.  For $m=2$, actually $M_2=K_4$ and two pebbles suffice, but
$C_2\cart K_2$ and $P_2\cart K_2$ degenerate to 4-cycles and need a third
pebble.

For the lower bound, we use induction on $m$.  For $m=1$ and $m=2$, note that
$\optpeb{P_m\cart K_2}=m+1$.  Now consider $m\ge 3$.  Since
$P_m\cart K_2\esub C_m\cart K_2$, it suffices to prove the lower bound for
$C_m\cart K_2$.  The argument for $C_m\cart K_2$ is valid also for $M_m$.

Consider an optimal solvable distribution $D$ with $|D|\le m$; we show that
equality holds.  If some pebbling sequence from $D$ results in a rung having
two pebbles, then collapsing that rung to a vertex yields a graph and
distribution under which the resulting vertex is $2$-reachable, so we say that
the rung is $2$-reachable under $D$.  If at least $m-1$ rungs are
$2$-reachable, then collapsing each rung to a vertex yields a distribution $D'$
on $C_m$ under which $m-1$ vertices are $2$-reachable.  Lemma~\ref{cycl2solv}
then yields $|D|=|D'|\ge m$.

Now suppose that at least two rungs $R$ and $R'$ are not $2$-reachable under
$D$.  The pebbles that arrive in pebbling sequences to reach the two vertices
of $R$ arrive from the same direction; otherwise, since no pebble can ever
emerge from $R'$, the two pebbling sequences can be performed independently
and $R$ is $2$-reachable.

Since both sequences reach $R$ from the same side, and no pebble can emerge
from $R$ to the other side (because $R$ is not $2$-reachable), $D$ remains
solvable on the graph obtained by deleting the edges from $R$ to that rung.
If there are two nonadjacent rungs that are not 2-solvable, then doing this for
those two rungs splits $D$ into solvable distributions on $P_i\cart K_2$ and
$P_{m-i}\cart K_2$, for some $i$ with $1\le i\le m-1$.  The induction
hypothesis applies to both subgraphs, and we obtain $|D|\ge m$.

In the remaining case, there are exactly two rungs $R$ and $R'$ that are not
$2$-reachable, and they are consecutive.  A rung that is not $2$-reachable is
unoccupied, because if there is one pebble on it, then the sequence to reach
the other vertex requires bringing another pebble to the rung.  Furthermore,
the pebbling sequences that move two pebbles to other rungs cannot use vertices
in $R$ or $R'$, since they are not 2-reachable.

Therefore, deleting $R$ and $R'$ and collapsing the remaining rungs yields
a 2-solvable distribution $D'$ on $P_{m-2}$.  If $|D'|\ge m$, we have the
desired result.  Otherwise, $D'$ is a minimal 2-solvable distribution on
$P_{m-2}$.  We use the description of all such distributions, obtained in
Theorem~\ref{path2}.

Let $S$ be the rung other than $R'$ that neighbors $R$; in the collapsed path,
$S$ is an endpoint.  Under $D'$, $S$ can receive two pebbles from its neighbor
if the prime segment ends $0$, or one pebble from its neighbor to join its
original pebble if the segment ends with 1, or no pebbles to join its two
original pebbles if the segment ends with 2.  In no case can $S$ receive a
third pebble.  Also, each case leaves no choice in the uncollapsed original
distribution $D$ about which vertex of the rung $S$ receives the extra pebble
or pair.  Without getting a third pebble to $S$ or being able to move two
pebbles to either vertex of $S$, it is not possible under $D$ to reach
each vertex of $R$.
\end{Proof}

\end{document}